\documentclass[8pt,a4paper]{amsart}
\usepackage[english]{babel}

\usepackage{amsmath,amssymb,amsthm}

\theoremstyle{plain}
\newtheorem{theorem}{Theorem}[section]

\theoremstyle{definition}
\newtheorem{remark}[theorem]{Remark}

\newcounter{equi1}

\begin{document}

\title
{Improving Three-point Iterative Methods for Solving Nonlinear Equations}

\author{Fayyaz Ahmad }
\author{D. García-Senz  }
\address{Dept. de Física i Enginyeria Nuclear, Universitat Politècnica de Catalunya, Comte d'Urgell 187, 08036 Barcelona, Spain } \email{fayyaz.ahmad@upc.edu, domingo.garcia@upc.edu}

\keywords{Non-linear equations, Iterative methods, Derivative free,  Steffensen's method,  Convergence order}



\begin{abstract}
  In this article, we report on sixth-order and seventh-order iterative methods for solving nonlinear equations. 
In particular sixth-order derivative-based and derivative-free iterative families  are constructed in such a way that they comprise a
wide class of sixth-order methods which were developed in the past years. Weighting functions are introduced to enhance the algorithmic 
efficiency whereas an appropriate parametric combination gives weight-age flexibility in between those weighting functions. The usage of 
weighting factors and weighting functions define a wide class of iterative schemes for solving nonlinear equations. The freedom to 
construct different parametric combinations as well as different forms of weighting functions makes the iterative schemes more accurate and 
flexible, It means that one can easily modify the scheme by changing weight functions and parametric combination. 
\end{abstract}

\maketitle

\section{Introduction}

Nonlinear algebraic equations are at the heart of many problems of nonlinear science. However, it is not always possible to find exact analytical solutions of 
these nonlinear equations which very often have to be solved numerically. Iterative methods are numerical algorithms which, starting from an initial guess in the neighborhood of a root of  the 
nonlinear equation, manage to refine it and  achieve convergence to the true root after several iterations  . We can divide them in two classes, namely  memory-based and without-memory 
iterative methods. From the point of view of efficiency and stability, multipoint iterative methods are superior than single-point methods. 
Kung and Traub conjectured \cite{13} that multipoint iterative methods without-memory have at most $2^n$ order of convergence for $n+1$ 
function evaluations per iteration. The further classification falls into derivative-free and non-derivative-free iterative methods. 
The multipoint iterative methods for nonlinear equations, which use derivatives, are efficient. The most famous derivative-based 
iterative method is Newton's method (\textbf{NM}) which is written as 
\begin{align} \label{eq1}
 x_{n+1} = x_n - \frac{f(x_n)}{f'(x_n)}, \ \ \ \  n=0,1,2,3,\cdots. 
\end{align}
It is quadratically convergent in some neighborhood of the root $\alpha$ of $f(x)=0$. The computational effort invested in the derivative of a function
is not always comparable to the evaluation of the function itself \cite{2,3,4,5,6,7,8}. Steffensen \cite{1} gave a quadratically convergent derivative-free iterative method
for nonlinear equations. The Steffensen's method (\textbf{SM}) is given by
\begin{align} \label{eq2}
\begin{cases}
  w_n     &= x_n-\kappa f(x_n), \\
  df(x_n) &= \frac{ f(x_n)-f(w_n) }{ \kappa f(x_n) } \\
  x_{n+1} &= x_n -  \frac{f(x_n)}{df(x_n)}.
\end{cases}
\end{align}
Clearly, $df(x_n)$ is a finite difference approximation of $f'(x_n)$. The \textbf{SM} approximation merely replaces the derivative evaluation by a function evaluation.
Actually both iterative method \textbf{NM} and \textbf{SM}  achieve the optimal order of convergence $2$, which is in agreement with the conjecture of 
Kung and Traub. The computational efficiency of an iterative method (IM) is frequently calculated by the Ostrowski-Traub's
\cite{14,15} efficiency index
\begin{align} \label{eqn1}
 E(IM) &= p^{1/d},
\end{align}
where $p$ is order of convergence and $d$ is function evaluations per iteration. The optimal computational efficiency is defined by expression 
\begin{align}\label{eqn2}
 E_{opt}^{(n)} &= 2^{n/(n+1)}.
\end{align}
We now briefly review a number of three step sixth-order convergent methods. Many three steps sixth-order
convergent methods use three function and one derivative evaluations. Sharma and Guha \cite{9} proposed the 
following sixth-order convergent family by introducing one parameter (\textbf{SG})
\begin{align} \label{eq3}
 \begin{cases}
  y_n &= x_n - \frac{f(x_n)}{f'(x_n)},\\
  z_n &= y_n - \frac{f(x_n)}{f(x_n)-2f(y_n)} \frac{f(y_n)}{f'(x_n)}, \\
  x_{n+1} &=  z_n- \frac{f(x_n)+a f(y_n)}{f(x_n)+(a-2) f(y_n)} \frac{f(z_n)}{f'(x_n)},
 \end{cases}
\end{align}
where $a\in \Re$. Neta \cite{10} uni-parametric sixth-order family consists of three step (\textbf{NT1})
\begin{align} \label{eq4}
 \begin{cases}
  y_n &= x_n - \frac{f(x_n)}{f'(x_n)},\\
  z_n &= y_n - \frac{f(x_n)+a f(y_n)}{f(x_n)+(a-2)f(y_n)} \frac{f(y_n)}{f'(x_n)},\\
  x_{n+1}  &= z_n - \frac{f(x_n)-f(y_n)}{f(x_n)-3f(y_n)} \frac{f(z_n)}{f'(x_n)}, 
 \end{cases}
\end{align}
where $a \in \Re$. In 2011, Neta \cite{11} provided a more efficient sixth-order scheme (\textbf{NT2})
\begin{align} \label{eq5}
 \begin{cases}
 y_n &= x_n - \frac{f(x_n)}{f'(x_n)}, \\
 z_n &= y_n - \frac{f(y_n)}{f'(x_n)} \frac{1}{\bigg[  1-\frac{f(y_n)}{f(x_n)}  \bigg]^2}, \\
 x_{n+1} &= z_n - \frac{f(z_n)}{f'(x_n)} \frac{1}{\bigg[  1-\frac{f(y_n)}{f(x_n)} -\frac{f(z_n)}{f(x_n)}    \bigg]^2}. 
 \end{cases}
\end{align}
Chun and Ham \cite{12} also presented a family of sixth-order iterative methods  (\textbf{CH})
\begin{align} \label{eq6}
 \begin{cases}
  y_n &= x_n - \frac{f(x_n)}{f'(x_n)}, \\
  z_n &= y_n -\frac{f(x_n)}{f(x_n)-2f(y_n)} \frac{f(y_n)}{f'(x_n)}, \\
  x_{n+1} &= z_n - H(\mu_n) \frac{f(z_n)}{f'(x_n)}, 
 \end{cases}
\end{align}
where real valued function $H(s)$ satisfies $H(0)=1$, $H'(0)=2$, and $\mu_n=f(y_n)/f(x_n)$. Grau et al. \cite{16} developed the sixth-order 
(\textbf{GR})
\begin{align} \label{eq7}
 \begin{cases}
  y_n &= x_n - \frac{f(x_n)}{f'(x_n)}, \\
  z_n &= y_n -\frac{f(x_n)}{f(x_n)-2f(y_n)} \frac{f(y_n)}{f'(x_n)}, \\
  x_{n+1} &= z_n - \frac{f(x_n)}{f(x_n)-2f(y_n)} \frac{f(z_n)}{f'(x_n)}.
 \end{cases}
\end{align}
Alicia cordero et al. \cite{17} constructed the following sixth-order scheme (\textbf{AL})
\begin{align} \label{eq8}
 \begin{cases}
  y_n &= x_n - \frac{f(x_n)}{f'(x_n)}, \\
  z_n &= x_n - \theta \frac{f(x_n)+f(y_n)}{f'(x_n)}-(1-\theta) \frac{f(x_n)}{f'(x_n)} \frac{f(x_n)}{f(x_n)-f(y_n)}, \\
  x_{n+1} &= z_n - \frac{f(z_n)}{f[z_n,y_n]+f[z_n,x_n,x_n](z_n-y_n)},
 \end{cases}
\end{align}
where $\theta \in \Re$.Finally, Sanjay K. Khattri et al. \cite{21} presented a paper, Unification of sixth-order iterative methods, which covers all the features  of 
sixth-order iterative methods (\textbf{SK1}) 
\begin{align} \label{eq13}
 \begin{cases}
  y_n &= x_n - \frac{f(x_n)}{f'(x_n)}, \\
  z_n &= y_n - \frac{f(y_n)}{f'(x_n)} \bigg(  1 + \sum\limits_{j=1}^m \bigg(  \frac{f(y_n)}{f(x_n)} \bigg)^j a_j  \bigg), \ \text{where}\ a_1=2,\\
  x_{n+1} &= z_n - \frac{f(z_n)}{f'(x_n)}  \bigg(1+\sum\limits_{k=1}^l \bigg(   \frac{\mu_1 f(y_n) + \mu_2 f(z_n)}{f(x_n)}   \bigg)^k  b_k \bigg), \ \text{where}\ b_1 = 2/\mu_1,
 \end{cases}
\end{align}
where, $\{a_j$, $b_k$, $\mu_1$, $\mu_2\}\in \Re$, $m \in \textbf{Z}^+ $, and $l\in \textbf{Z}^+$ are independent parameters.
Sanjay K. Khattri et al. \cite{18} also developed four parameter sixth-order derivative free family (\textbf{SK2})
\begin{align} \label{eq9}
 \begin{cases}
  y_n &= x_n -\kappa \frac{f(x_n)^2}{f(x_n)-f(x_n-\kappa f(x_n))},\\
  z_n &= y_n - \kappa \frac{f(x_n)f(y_n)}{f(x_n)-f(x_n-\kappa f(x_n))} \bigg[ 1+\frac{f(y_n)}{f(x_n)} + \alpha \bigg(\frac{f(y_n)}{f(x_n)} \bigg)^2 
  +\frac{f(y_n)}{f(x_n-\kappa f(x_n)) } + \beta \bigg(\frac{f(y_n)}{f(x_n-\kappa f(x_n)) }\bigg)^2\bigg], \\
  x_{n+1} &= z_n - \kappa \frac{f(x_n)f(z_n)}{f(x_n)-f(x_n-\kappa f(x_n))}\bigg[ 1+\frac{f(y_n)}{f(x_n)} + \alpha \bigg(\frac{f(y_n)}{f(x_n)} \bigg)^2 
  +\frac{f(y_n)}{f(x_n-\kappa f(x_n)) } + \beta \bigg(\frac{f(y_n)}{f(x_n-\kappa f(x_n)) }\bigg)^2 + \eta \frac{f(z_n)}{f(y_n)} \bigg] .
 \end{cases}
\end{align}
In \cite{19}, R. Thukral constructed Steffensen-Six.one method (\textbf{TS1}) and Steffensen-six.two method  (\textbf{TS2}), which are respectively
\begin{align} 
 \begin{cases}\label{eq10}
  w_n &= x_n + f(x_n), \\
  y_n &= x_n - \bigg( \frac{f(x_n)^2}{f(w_n)-f(x_n)}   \bigg), \\
  z_n &= y_n - \bigg(  1- \lambda^{-1}  \frac{f(y_n)}{f(w_n)} \bigg)^{-\lambda} \bigg( \frac{x_n-y_n}{f(x_n)-f(y_n)}  \bigg) f(y_n), \\
  x_{n+1} &= z_n - \bigg(  1- \lambda^{-1}  \frac{f(y_n)}{f(w_n)} \bigg)^{-\lambda} \bigg( \frac{x_n-y_n}{f(x_n)-f(y_n)}  \bigg) f(z_n),
 \end{cases}
 \intertext{and}
 \begin{cases} \label{eq11}
    w_n &= x_n + f(x_n), \\
  y_n &= x_n - \bigg( \frac{f(x_n)^2}{f(w_n)-f(x_n)}   \bigg), \\
  z_n &= y_n - \bigg(  1- \lambda^{-1}  \frac{f(y_n)}{f(w_n)} \bigg)^{-\lambda} \bigg( \frac{x_n-y_n}{f(x_n)-f(y_n)}  \bigg) f(y_n), \\
  x_{n+1} &= z_n - \bigg( \frac{z_n-y_n}{f(z_n)-f(y_n)}  \bigg) f(z_n),
 \end{cases}
\end{align}
where $\lambda \in \Re \slash \{ 0 \} $. F. Soleymani et al. \cite{20} discussed the following sixth-order derivative-free method (\textbf{FS1})
\begin{align} \label{eq12}
 \begin{cases}
  w_n &= x_n + f(x_n), \\
  y_n &= x_n - \frac{f(x_n)}{f[x_n,w_n]}, \\
  z_n &= x_n - \frac{f(x_n)}{f[x_n,w_n]} \bigg[ 1+\frac{f(y_n)}{f(x_n)}\bigg( 1+2 \frac{f(y_n)}{f(x_n)}  \bigg)  \bigg],\\
  x_{n+1} &= z_n - \frac{f(z_n)}{f[y_n,z_n]} \bigg( 1 - \frac{1+f[x_n,w_n]}{f[x_n,w_n]}  \frac{f(z_n)}{f(w_n)}\bigg),
 \end{cases}
\end{align}
where $f[y_n,z_n] = \frac{f(z_n)-f(y_n)}{z_n-y_n}$, $f[x_n,w_n] = \frac{f(w_n)-f(x_n)}{w_n-x_n}$. F. Soleymani \cite{22} also presented an other 
sixth-order derivative-free scheme (\textbf{FS2})
\begin{align} \label{eq14}
\begin{cases}
 w_n &= x_n - \beta f(x_n),\\
 y_n &= x_n - \frac{f(x_n)}{f[x_n,w_n]}, \\
 z_n &= y_n - \frac{f(y_n)}{f[w_n,y_n]},\\
 x_{n+1} &= z_n - \frac{f(z_n)}{f[w_n,z_n]+f[z_n,y_n]-f[w_n,y_n]},
\end{cases}
\end{align}
where $\beta \in \Re \slash \{ 0 \}$.In 2013, F. Soleymani \cite{26} developed some efficient seventh-order derivative-free methods.
First is (\textbf{FS3-1})
\begin{align}\label{eq59}
 \begin{cases}
  w_n &= x_n + f(x_n), \\
  y_n &=x_n-\frac{f(x_n)}{f[x_n,w_n]}, \\
  z_n &=y_n-\frac{f(y_n)}{f[x_n,y_n]+f[y_n,w_n]-f[x_n,w_n]},\\
  x_n &= z_n - \frac{f(z_n)}{f[x_n,z_n]}\Bigg[ 1+\frac{f(y_n)}{f(w_n)}+\frac{f(z_n)}{f(y_n)}+\Bigg\{  \frac{2+f[x_n,w_n]}{(1+f[x_n,w_n])^2}\Bigg\} 
  \Bigg(\frac{f(y_n)}{f(x_n)}\Bigg)^2 + \gamma \frac{f(z_n)}{f(x_n)} + \delta \frac{f(z_n)}{f(w_n)} \Bigg],
 \end{cases}
\end{align}
where $\gamma$ and $\delta$ are real numbers. Second is (\textbf{FS3-2})
\begin{align}\label{eq60}
 \begin{cases}
  w_n &= x_n - f(x_n), \\
  y_n &=x_n-\frac{f(x_n)}{f[x_n,w_n]}, \\
  z_n &=y_n-\frac{f(y_n)}{f[x_n,y_n]+f[y_n,w_n]-f[x_n,w_n]},\\
  x_n &= z_n - \frac{f(z_n)}{f[x_n,z_n]}\Bigg[ 1+\frac{f(y_n)}{f(w_n)}+\frac{f(z_n)}{f(y_n)}+\Bigg\{  \frac{2-f[x_n,w_n]}{(1-f[x_n,w_n])^2}\Bigg\} 
  \Bigg(\frac{f(y_n)}{f(x_n)}\Bigg)^2 + \rho \frac{f(z_n)}{f(x_n)} + \tau \frac{f(z_n)}{f(w_n)} \Bigg],
 \end{cases}
\end{align}
where $\rho$ and $\tau$ are real numbers. Third is (\textbf{FS4-1})
\begin{align}\label{eq61}
 \begin{cases}
  w_n &= x_n + f(x_n), \\
  y_n &=x_n-\frac{f(x_n)}{f[x_n,w_n]}, \\
  z_n &=y_n-\frac{f(y_n)}{f[x_n,y_n]+f[y_n,w_n]-f[x_n,w_n]},\\
  x_n &= z_n - \frac{f(z_n)}{f[w_n,z_n]}\Bigg[ 1+\frac{f(z_n)}{f(y_n)}+\frac{f(y_n)}{f(x_n)}+\Big\{ 2+f[x_n,w_n](3+f[x_n,w_n]) \Big\} 
  \Bigg(\frac{f(y_n)}{f(w_n)}\Bigg)^2 + \omega \frac{f(z_n)}{f(x_n)} + \phi \frac{f(z_n)}{f(w_n)} \Bigg],
 \end{cases}
\end{align}
where $\omega$ and $\phi$ are real numbers. Fourth is (\textbf{FS4-2})
\begin{align}\label{eq62}
 \begin{cases}
  w_n &= x_n - f(x_n), \\
  y_n &=x_n-\frac{f(x_n)}{f[x_n,w_n]}, \\
  z_n &=y_n-\frac{f(y_n)}{f[x_n,y_n]+f[y_n,w_n]-f[x_n,w_n]},\\
  x_n &= z_n - \frac{f(z_n)}{f[w_n,z_n]}\Bigg[ 1+\frac{f(z_n)}{f(y_n)}+\frac{f(y_n)}{f(x_n)}+\Big\{ 2+f[x_n,w_n](-3+f[x_n,w_n]) \Big\} 
  \Bigg(\frac{f(y_n)}{f(w_n)}\Bigg)^2  \Bigg].
 \end{cases}
\end{align}
\section{A More general formulation of sixth-order iterative methods}
In this section, we go beyond the iterative method (\textbf{SK1}) and 
present a more general formulation which covers most of the sixth-order iterative methods based on one
derivative evaluation $f'(x)$. The proposed iterative scheme is (\textbf{FD1})
\begin{align}\label{eq15}
 \begin{cases}
  y_n &= x_n - \frac{f(x_n)}{f'(x_n)},\\
  z_n &= y_n-G (t_1) \frac{f(y_n)}{f'(x_n)},\\
  x_{n+1} &= z_n - H (t_1,t_2,t_3) \frac{f(z_n)}{f'(x_n)},
 \end{cases}
\end{align}
where $t_1=\frac{f(y_n)}{f(x_n)}$, $t_2=\frac{f(z_n)}{f(x_n)}$ and $t_3=\frac{f(z_n)}{f(y_n)}$. we state the following theorem.
\newtheorem{mydef1}{Theorem}
\begin{mydef1}
 Let  $f: D\subseteq \Re \rightarrow \Re $ be a sufficiently differentiable function, and $\alpha \in D$ is a simple root of $f(x)=0$, 
 for an open interval $D$. If $x_0$ is chosen sufficiently close to $\alpha$, then iterative scheme (\ref{eq15}) converges to $\alpha$. If $G$ 
 and $H$ satisfy 
  \begin{align} \label{eq16}
  G(0)&=1,\ dG /dt_1 (0) =2,\ 
  H(0, 0, 0)=1,\ \partial H / \partial t_1 (0, 0, 0) =2,   
 \end{align}
 and both $G$, $H$ have  bounded higher order derivatives then the iterative scheme (\ref{eq15}) has convergence order at-least six and 
 along with  conditions (\ref{eq16}) if
 \begin{align} \label{eq27}
 \partial^2 H / \partial t_1^2 (0, 0, 0) =d^2 G /dt_1^2 (0)+ 2,\   \partial H / \partial t_3 (0, 0, 0) =1  
 \end{align}
 then the iterative scheme (\ref{eq15}) has convergence order at least seven.
\end{mydef1}
\begin{proof}
 let  $\alpha$ be the simple root of $f(x)=0$ and  $e_n = x_n-\alpha$. Further we denote $c_1 = f'(\alpha)$ and $c_k = 
 \frac{e_n^k}{k!}  \frac{f^{(k)}(\alpha)}{f'(\alpha)}$, where $k=2, 3, \cdots$. The Tyalor's expression of $f$ and $f'$ around $\alpha$ in terms of $e_n$:
 \begin{align} 
 f(x_n)&= c_1 [e_n + c_2 e_n^2 + c_3 e_n^3 + c_4 e_n^4+ c_5 e_n^5 +c_6 e_n^6 +O(e_n^7)]  \label{eq21} \\
 f'(x_n)&= c_1 [1 + 2 c_2 e_n + 3 c_3 e_n^2 + 4 c_4 e_n^3 +5 c_5 e_n^4 +6 c_6 e_n^5 + O(e_n^6)]  \label{eq22} 
 \end{align}
 Substituting (\ref{eq21}) and (\ref{eq22}) into $y_n$, we obtain 
 \begin{align}
  y_n &= x_n - \frac{f(x_n)}{f'(x_n)},\  
  y_n-\alpha = x_n-\alpha  - \frac{f(x_n)}{f'(x_n)} = e_n - \frac{f(x_n)}{f'(x_n)} \nonumber \\
  &=    c_2 e_n^2+ (2 c_3-2 c_2^2 ) e_n^3 +(3 c_4-7 c_2 c_3+4 c_2^3 ) e_n^4  +( 4 c_5-10 c_2 c_4-6 c_3^2+20 c_3 c_2^2-8 c_2^4 ) e_n^5 \nonumber \\
  &+(-17 c_3 c_4+33 c_2 c_3^2-52 c_3 c_2^3+28 c_4 c_2^2-13 c_2 c_5+5 c_6+16 c_2^5 ) e_n^6+O(e_n^7)    \label{eq23}  
 \end{align}
  Expansion of $f(y_n)$ around $\alpha$ gives
  \begin{align}
   f(y_n) &=  c_1[  c_2 e_n^2 - 2  (-c_3+c_2^2)  e_n^3 + (3 c_4-7 c_2 c_3+5 c_2^3)  e_n^4  -2  (-2 c_5+5 c_2 c_4+3 c_3^2-12 c_3 c_2^2+6 c_2^4)  e_n^5 \nonumber \\
   &  c_1 (37  c_2  c_3^2-73  c_3  c_2^3+28  c_2^5+34  c_4  c_2^2-17  c_3  c_4-13  c_2  c_5+5  c_6) e_n^6+O(e_n^7)]   \label{eq24}
  \end{align}
   We define $G(t_1) = 1 + 2 t_1 + M_0 t_1^2 + M_1 t_1^3 + M_2 t_1^4+ \text{higher powers in $t_1$}$ 
   , where $M_0$, $M_1$ and $M_2$ are finite real numbers. By using (\ref{eq21})-(\ref{eq24}), the expressions for $z_n$ and $f(z_n)$
   are 
   \begin{align}
    z_n &= y_n -G(t_1) \frac{f(y_n)}{f'(x_n)} \nonumber \\
    z_n-\alpha &=  -c_2 (-5 c_2^2+M_0 c_2^2+c_3)  e_n^4 + (-2 c_2 c_4-2 c_3^2+32 c_3 c_2^2-36 c_2^4-M_1 c_2^4-6 M_0 c_2^2 c_3+10 M_0 c_2^4) e_n^5 \nonumber \\
    & (-12 M_0  c_2  c_3^2-9 M_0  c_2^2  c_4+74 M_0  c_2^3  c_3-62 M_0  c_2^5-7  c_3  c_4+66  c_2  c_3^2-262  c_3  c_2^3+48  c_4  c_2^2-3  c_2  c_5\nonumber \\
    & +170  c_2^5-M_2  c_2^5+13 M_1  c_2^5-8 M_1  c_2^3  c_3) e_n^6+O(e_n^7) \label{eq25} \\
    f(z_n )    &=  - c_1  c_2 (-5  c_2^2+M_0  c_2^2+ c_3)  e_n^4 +  c_1 (-2  c_2  c_4-2  c_3^2+32  c_3  c_2^2-36  c_2^4-M_1  c_2^4-6 M_0  c_2^2  c_3 \nonumber \\
    &+10 M_0  c_2^4)  e_n^5 +c_1(-12 M_0  c_2  c_3^2-9 M_0  c_2^2  c_4+74 M_0  c_2^3  c_3-62 M_0  c_2^5-7  c_3  c_4+66  c_2  c_3^2-262  c_3  c_2^3\nonumber \\
    & +48  c_4  c_2^2-3  c_2  c_5+170  c_2^5-M_2  c_2^5+13 M_1  c_2^5-8 M_1  c_2^3  c_3)e_n^6+O(e_n^7) \label{eq26} 
   \end{align}
   Now, we define $H(t_1,t_2 ,t_3) =1+ 2 t_1 +  R_0 t_2 + R_1 t_3 + R_2  t_1^2 +  R_3  t_2^2+ R_4  t_3^2+ \text{higher powers in $t_1$, $t_2$, $t_3$}$. 
    Again, substituting (\ref{eq22}), (\ref{eq25}) and (\ref{eq26}) in (\ref{eq15}), we obtain
    \begin{align} \label{eq28}
     x_{n+1}-\alpha &= z_n -\alpha - H(t_1,t_2,t_3) \frac{f(y_n)}{f'(x_n)} \nonumber\\
     & =  - c_2 (-5  c_2^2+M_0  c_2^2+ c_3) (R_1 M_0  c_2^2-5 R_1  c_2^2+6  c_2^2-R_2  c_2^2+R_1  c_3- c_3) e_n^6  \nonumber \\
     & +( 10 R_0 c_2^6 M_0-R_0 c_2^6 M_0^2+10 R_1 M_1 c_2^6-172 R_1 c_2^6 M_0+20 R_1 M_0^2 c_2^6+2 M_0 c_2^3 c_4+6 M_0 c_2^2 c_3^2\nonumber \\
     & -68 M_0 c_2^4 c_3+88 M_0 c_2^6+124 R_1 M_0 c_2^4 c_3-12 R_1 M_0 c_2^2 c_3^2-10 R_1 M_0^2 c_2^4 c_3-2 R_1 M_1 c_2^6 M_0 -2 R_1 M_1 c_2^4 c_3\nonumber \\
     & -2 R_0 c_2^4 M_0 c_3-4 R_1 c_4 M_0 c_2^3-22 c_4 c_2^3-66 c_3^2 c_2^2+366 c_3 c_2^4+2 c_3^3-356 c_2^6+4 c_4 c_2 c_3-6 M_1 c_2^6\nonumber \\
     & -25 R_6 c_2^6+10 R_6 c_2^6 M_0+10 R_6 c_2^4 c_3-R_6 c_2^6 M_0^2-R_6 c_2^2 c_3^2+M_1 c_2^4 c_3+76 R_2 c_2^6-2 R_6 c_2^4 M_0 c_3\nonumber \\
     & -4 c_2 R_1 c_4 c_3-342 R_1 c_3 c_2^4+20 R_1 c_4 c_2^3+10 R_0 c_2^4 c_3+64 c_2^2 R_1 c_3^2-R_0 c_2^2 c_3^2-25 R_0 c_2^6+360 R_1 c_2^6\nonumber \\
     & -2 R_1 c_3^3-60 R_2 c_2^4 c_3+6 R_2 c_2^2 c_3^2+2 R_2 c_2^3 c_4-18 R_2 c_2^6 M_0+R_2 c_2^6 M_1+10 R_2 c_2^4 M_0 c_3) e_n^7 +O(e_n^8)
    \end{align}
By equating  one of the factor of coefficient of $ e_n^6 $ to zero, we have 
\begin{align} 
 R_1 M_0  c_2^2-5 R_1  c_2^2+6  c_2^2-R_2  c_2^2+R_1  c_3- c_3 =0. \label{eq29} 
\end{align}
Now by equating the coefficient of $c_3$ and $c_2^2$ in (\ref{eq29}), we find
\begin{align} \label{eq30}
\begin{cases}
 &R_1 -1 = 0,  \  \text{(coefficient of $c_3$) }\\
 &R_1 M_0  -5 R_1  +6  -R_2 =0, \   \text{ (coefficient of $c_2^2$) }
\end{cases}
\end{align}
where $M_0=1/2\ d^2 G/dt_1^2 (0) $, $R_1 =  \partial H / \partial t_3 (0, 0, 0) $ and $ R_2 = 1/2\ \partial^2 H / \partial t_1^2 (0, 0, 0) $. Which completes
the proof.
\end{proof}
\begin{remark}
 It is obvious that sixth-order family (\ref{eq15}) can construct all the presented iterative methods in (\ref{eq3})-(\ref{eq7}) and (\ref{eq13}).
\end{remark}
(\ref{eq8}) can be written as (\textbf{FD2})
\begin{align}\label{eq17}
 \begin{cases}
  y_n &= x_n - \frac{f(x_n)}{f'(x_n)}, \\
  z_n &= x_n - A(t_1) \frac{f(y_n)}{f'(x_n)},\\
  x_{n+1} &= z_n - \frac{f(z_n)}{f[z_n,y_n]+f[z_n,x_n,x_n](z_n-y_n)},
 \end{cases}
\end{align}
where $t_1=\frac{f(y_n)}{f(x_n)}$ and $A(t_1)=\bigg(  \frac{1-\theta t_1}{1-t_1} \bigg)$,  with $\theta \in \Re$ and $f[z_n,x_n,x_n] = \frac{f[z_n,x_n]-f'(x_n)}{z_n-x_n}$. The more general form  of (\ref{eq17})
is stated in the following theorem.
\newtheorem{mydef2}[mydef1]{Theorem}
\begin{mydef2}
 Let  $f: D\subseteq \Re \rightarrow \Re $ be a sufficiently differentiable function, and $\alpha \in D$ is a simple root of $f(x)=0$, 
 for an open interval $D$. If $x_0$ is chosen sufficiently close to $\alpha$, then iterative scheme (\ref{eq17}) converges to $\alpha$. If $A$
  satisfies
  \begin{align} \label{eq18}
  A(0)&=1
 \end{align}
 then the iterative scheme (\ref{eq17}) has convergence order at least six and if along with (\ref{eq18}), it satisfies
 \begin{align}\label{eq31}
   dA/dt_1 (0) =2 
 \end{align}
  then  (\ref{eq17}) has convergence order at least seven.
\end{mydef2}
\begin{proof}
 We define $A(t_1) = 1+M_0 t_1 +M_1 t_1^2 +M_2 t_1^3 + \text{higher powers in $s$}$. By using the (\ref{eq21})-(\ref{eq24}), we obtain
 \begin{align}  
  z_n-\alpha &= y_n-\alpha -(1+M_0 t_1 +M_1 t_1^2 +M_2 t_1^3) \frac{f(y_n)}{f'(x_n)}  \nonumber \\ 
  &= (-c_2^2 (-2+M_0)) e_n^3+ c_2 (-9 c_2^2+7 M_0 c_2^2-M_1 c_2^2+7 c_3-4 M_0 c_3)  e_n^4  +(10 c_2 c_4+6 c_3^2-44 c_3 c_2^2\nonumber \\
  & +30 c_2^4-6 M_0 c_2 c_4+38 M_0 c_2^2 c_3-33 M_0 c_2^4-M_2 c_2^4-6 M_1 c_2^2 c_3+10 M_1 c_2^4-4 M_0 c_3^2) e_n^5 + (17  c_3  c_4\nonumber \\
  & -70  c_2  c_3^2+188  c_3  c_2^3-62  c_4  c_2^2+13  c_2  c_5-88  c_2^5+74 M_1  c_2^3  c_3-9 M_1  c_2^2  c_4-12 M_1  c_2  c_3^2-62 M_1  c_2^5\nonumber \\
  & +13 M_2  c_2^5-8 M_2  c_2^3  c_3+129 M_0  c_2^5-12 M_0  c_3  c_4+68 M_0  c_2  c_3^2-225 M_0  c_3  c_2^3+55 M_0  c_4  c_2^2-8 M_0  c_2  c_5) e_n^6 \nonumber \\
  & +O(e_n^7), \label{eq32} \\
  f(z_n) &= c_1 [(-c_2^2 (-2+M_0)) e_n^3+ c_2 (-9 c_2^2+7 M_0 c_2^2-M_1 c_2^2+7 c_3-4 M_0 c_3)  e_n^4  +(10 c_2 c_4+6 c_3^2-44 c_3 c_2^2\nonumber \\
  & +30 c_2^4-6 M_0 c_2 c_4+38 M_0 c_2^2 c_3-33 M_0 c_2^4-M_2 c_2^4-6 M_1 c_2^2 c_3+10 M_1 c_2^4-4 M_0 c_3^2) e_n^5 + (17  c_3  c_4\nonumber \\
  & -70  c_2  c_3^2+188  c_3  c_2^3-62  c_4  c_2^2+13  c_2  c_5-88  c_2^5+74 M_1  c_2^3  c_3-9 M_1  c_2^2  c_4-12 M_1  c_2  c_3^2-62 M_1  c_2^5\nonumber \\
  & +13 M_2  c_2^5-8 M_2  c_2^3  c_3+129 M_0  c_2^5-12 M_0  c_3  c_4+68 M_0  c_2  c_3^2-225 M_0  c_3  c_2^3+55 M_0  c_4  c_2^2-8 M_0  c_2  c_5) e_n^6 \nonumber \\
  & +O(e_n^7)]. \label{eq33}  
 \end{align}
  By substituting values from (\ref{eq21})-(\ref{eq24}), (\ref{eq32}) and (\ref{eq33}) in $x_{n+1}$ from (\ref{eq17}), we get 
  \begin{align} \label{eq34}
  x_{n+1}-\alpha  &= z_n -\alpha  - \frac{f(z_n)}{f[z_n,y_n]+f[z_n,x_n,x_n](z_n-y_n)} \nonumber \\
  e_{n+1}  &=   c_2^3 (-2+M_0) (M_0  c_2^2+2  c_3-2  c_2^2) e_n^6  -c_2^2 (36 c_2^4-46 M_0 c_2^4+4 M_1 c_2^4-2 M_1 c_2^4 M_0+14 M_0^2 c_2^4\nonumber\\
  & -10 M_0^2 c_3 c_2^2-2 M_1 c_2^2 c_3-64 c_3 c_2^2+57 M_0 c_2^2 c_3+6 c_2 c_4-3 M_0 c_2 c_4+22 c_3^2-12 M_0 c_3^2) e_n^7 + O(e_n^8).   
  \end{align}
  Clearly, if $ dA(t_1)/dt_1(0)= M_0 = 2$ then convergence order  is seven otherwise is six.
\end{proof}
Now we present derivative-free optimal fourth-order iterative schemes (\textbf{FD3}): 
\begin{align} \label{eq19}
 \begin{cases}
  w_n     & = x_n-\kappa f(x_n), \\
  y_n     & = x_n - \frac{f(x_n)}{f[x_n,w_n]},\\
  x_{n+1}     & = y_n -  \Big[  \frac{g_0}{f[y_n,w_n]} G_0( t_1)   +    \frac{g_1}{f[y_n,x_n]} 
              G_1( t_2) +  \frac{g_2}{f[x_n,w_n]} G_2(t_1,t_2)   \Big],
 \end{cases}
\end{align}
where $t_1=\frac{f(y_n)}{f(x_n)}$, $t_2 = \frac{f(y_n)}{f(w_n)}$.
The convergence of (\ref{eq19}) is given in the following theorem.
\newtheorem{mydef3}[mydef1]{Theorem}
\begin{mydef3}
 Let  $f: D\subseteq \Re \rightarrow \Re $ be a sufficiently differentiable function, and $\alpha \in D$ is a simple root of $f(x)=0$, 
 for an open interval $D$. If $x_0$ is chosen sufficiently close to $\alpha$, then iterative scheme (\ref{eq19}) converges to $\alpha$. If
 $G_0(t_1)$, $G_1(t_2)$, $G_2(t_1,t_2)$, $g_0$, $g_1$ and $g_2$ satisfy 
  \begin{align} \label{eq20}
  \begin{cases}
  &G_0(0))=1,\  G_1(0)=1,   \ G_2(0,0)=1, \\   
  & dG_0/dt_1(0)=1,\ dG_1/dt_2(0)=1,\ \partial G_2 / \partial t_1 (0,0) =1,\ \partial G / \partial t_2 (0,0) =1,  \\
  & g_0+g_1+g_2=1 
  \end{cases}
 \end{align}
  and all higher order derivatives of $G_i,\ i=0,1,2$ are bounded then the iterative scheme (\ref{eq19}) has convergence order at least four.
\end{mydef3}
\begin{proof}
 Error term at nth-step is denoted by $e_n = x_n-\alpha$. we call $c_1 = f'(\alpha)$ and $c_k = 
 \frac{e_n^k}{k!}  \frac{f^{(k)}(\alpha)}{f'(\alpha)}$, where $k=2, 3, \cdots$. Taylor's expressions give:
 \begin{align} 
  f(x_n)&= c_1 [e_n + c_2 e_n^2 + c_3 e_n^3 + c_4 e_n^4 + c_5 e_n^5 + c_6 e_n^6 + O(e_n^7)]. \label{eq35} \\
  f(w_n) &= -c_1 (-1+ \kappa c_1) e_n + c_1 c_2 (-3  \kappa c_1+1+ \kappa^2 c_1^2) e_n^2  -c_1 (4  \kappa c_1 c_3-c_3-3 c_3  \kappa^2 c_1^2+c_3  \kappa^3 c_1^3+2 c_2^2  \kappa c_1 \nonumber \\
  &-2  \kappa^2 c_1^2 c_2^2) e_n^3+  c_1 (-5 \kappa  c_1  c_4-5  c_2 \kappa  c_1  c_3+8 \kappa^2  c_1^2  c_2  c_3+\kappa^2  c_1^2  c_2^3+ c_4+6  c_4 \kappa^2  c_1^2-4  c_4 \kappa^3  c_1^3+ c_4 \kappa^4  c_1^4\nonumber \\
  & -3  c_2 \kappa^3  c_1^3  c_3) e_n^4 +\cdots +O(e_n^7) \label{eq36}.
  \end{align}
  By using (\ref{eq35}) and (\ref{eq36}), we get
 \begin{align}
  f[x_n,w_n]& =\frac{f(x_n)-f(w_n)}{x_n-w_n} = -c_1 c_2 (-2+\kappa c_1)  e_n +  c_1 (3  c_3-3 \kappa  c_1  c_3+ c_3 \kappa^2  c_1^2- c_2^2 \kappa  c_1) e_n^2 - c_1 (-4  c_4\nonumber \\
  & +4  c_2 \kappa  c_1  c_3+6 \kappa  c_1  c_4-4  c_4 \kappa^2  c_1^2+ c_4 \kappa^3  c_1^3-2 \kappa^2  c_1^2  c_2  c_3) e_n^3 +  c_1 (2 \kappa ^2  c_1^2  c_3^2+8 \kappa ^2  c_1^2  c_4  c_2 \nonumber \\
  &-10 \kappa   c_1 c_5+5  c_5-7  c_4 \kappa   c_1  c_2-3 \kappa   c_1  c_3^2+10  c_5 \kappa ^2  c_1^2-5  c_5 \kappa ^3  c_1^3+ c_5 \kappa ^4  c_1^4-3  c_4 \kappa ^3  c_1^3  c_2+\kappa ^2  c_1^2  c_2^2  c_3) e_n^4 \nonumber \\
  & +\cdots +O(e_n^7) \label{eq37}
   \end{align}
  By substituting (\ref{eq35}) and (\ref{eq37}) in $y=x-f(x_n)/f[x_n,w_n]$, we obtain
  \begin{align}\label{eq38}
   y_n-\alpha &= -c_2 (-1+\kappa c_1)  e_n^2 +( 2  c_3-3 \kappa  c_1  c_3+ c_3 \kappa^2  c_1^2+2  c_2^2 \kappa  c_1-2  c_2^2-\kappa^2  c_1^2  c_2^2 ) e_n^3 + (3 c_ 4+10 c_ 2 \kappa c_ 1 c_ 3-6 \kappa c_ 1 c_ 4\nonumber \\
   & +4 c_ 4 \kappa^2 c_ 1^2-c_ 4 \kappa^3 c_ 1^3-7 \kappa^2 c_ 1^2 c_ 2 c_ 3-7 c_ 2 c_ 3-5 \kappa c_ 1 c_ 2^3+2 c_ 2 \kappa^3 c_ 1^3 c_ 3+3 \kappa^2 c_ 1^2 c_ 2^3+4 c_ 2^3-\kappa^3 c_ 1^3 c_ 2^3) e_n^4 +\cdots+O(e_n^7)
  \end{align}
   Taylor's Expansion of $f$ around $y_n$ by using (\ref{eq38}) is 
   \begin{align} \label{eq39}
    f(y_n) &= c_1( -c_2 (-1+\kappa c_1)  e_n^2 +( 2  c_3-3 \kappa  c_1  c_3+ c_3 \kappa^2  c_1^2+2  c_2^2 \kappa  c_1-2  c_2^2-\kappa^2  c_1^2  c_2^2 ) e_n^3 + (3 c_ 4+10 c_ 2 \kappa c_ 1 c_ 3-6 \kappa c_ 1 c_ 4\nonumber \\
   & +4 c_ 4 \kappa^2 c_ 1^2-c_ 4 \kappa^3 c_ 1^3-7 \kappa^2 c_ 1^2 c_ 2 c_ 3-7 c_ 2 c_ 3-5 \kappa c_ 1 c_ 2^3+2 c_ 2 \kappa^3 c_ 1^3 c_ 3+3 \kappa^2 c_ 1^2 c_ 2^3+4 c_ 2^3-\kappa^3 c_ 1^3 c_ 2^3) e_n^4 +\cdots+O(e_n^7))
   \end{align}
   Now, we define 
   \begin{align}
    \begin{cases} \label{eq40}
     G_0 (t_1)     &= 1+  t_1 + L_1 t_1^2+ \text{higher powers in $t_1$}, \\
     G_1 (t_2)     &= 1+ t_2 + M_1 t_2^2+ \text{higher powers in $t_2$}, \\
     G_2 (t_1,t_2) &= 1+ t_1 + t_2 + N_1 t_1^2   + N_2 t_2^2+ \text{higher powers in $t_1$, $t_2$},\\
    \end{cases}
   \end{align}
    where $L_1$, $M_1$, $N_1$ and $N_2$ are real numbers. By substituting (\ref{eq35})-(\ref{eq39}) in $x_{n+1}= y_n -\Big[ \frac{g_0}{f[y_n,w_n]} G_0(t_1) + \frac{g_1}{f[y_n,x_n]}  G_1(t_2)
             + \frac{g_2}{f[x_n,w_n]}  G_2 (t_1,t_2)  \Big]$, we obtain
    \begin{align}\label{eq41}
     e_{n+1} &=  (-1+g_0+g_1+g_2) c_ 2 (-1+ \kappa c_ 1) e_n^2  +\nonumber \\
     & (-1+g_0+g_1+g_2) (-2 c_3+3  \kappa c_1 c_3-c_3  \kappa^2 c_1^2-2 c_2^2  \kappa c_1+2 c_2^2 + \kappa^2 c_1^2 c_2^2)  e_n^3 + \nonumber \\
     & (-3 g_0 c_4-g_0 c_2^3+6 g_0 c_2 c_3-7 \kappa^2 c_1^2 c_2 c_3-6 \kappa c_1 c_4+3 c_4+3 g_2 N_1 c_2^3 \kappa c_1-3 g_2 N_1 c_2^3 \kappa^2 c_1^2\nonumber \\
     & +g_2 N_1 c_2^3 \kappa^3 c_1^3+g_2 N_2 c_2^3 \kappa c_1+g_1 M_1 c_2^3 \kappa c_1+g_0 \kappa^3 c_1^3 L_1 c_2^3+3 g_0 L_1 c_2^3 \kappa c_1-3 g_0 L_1 c_2^3 \kappa^2 c_1^2+3 \kappa^2 c_1^2 c_2^3\nonumber \\
     & +4 c_4 \kappa^2 c_1^2-c_4 \kappa^3 c_1^3-5 \kappa c_1 c_2^3-\kappa^3 c_1^3 c_2^3+10 c_2 \kappa c_1 c_3+2 c_2 \kappa^3 c_1^3 c_3-7 c_2 c_3+4 c_2^3-5 g_2 \kappa c_1 c_2^3+g_1 \kappa^3 c_1^3 c_2^3\nonumber \\
     & +6 g_0 \kappa c_1 c_4-4 g_0 c_4 \kappa^2 c_1^2+g_0 c_4 \kappa^3 c_1^3-2 g_0 \kappa c_1 c_2^3+2 g_0 \kappa^2 c_1^2 c_2^3-g_2 N_2 c_2^3-g_1 M_1 c_2^3-g_2 N_1 c_2^3-g_0 L_1 c_2^3\nonumber \\
     & -8 g_0 c_2 \kappa c_1 c_3+6 g_0 \kappa^2 c_1^2 c_2 c_3-2 g_0 c_2 \kappa^3 c_1^3 c_3+6 g_1 \kappa c_1 c_4+6 g_1 \kappa^2 c_1^2 c_2 c_3-8 g_1 c_2 \kappa c_1 c_3-2 g_1 c_2 \kappa^3 c_1^3 c_3\nonumber \\
     & +6 g_2 \kappa^2 c_1^2 c_2 c_3-8 g_2 c_2 \kappa c_1 c_3-2 g_2 c_2 \kappa^3 c_1^3 c_3+6 g_1 c_2 c_3+6 g_2 c_2 c_3-3 g_1 c_4-g_1 c_2^3-3 g_2 c_4+g_2 c_2^3\nonumber \\
     & -g_1 \kappa^2 c_1^2 c_2^3-4 g_1 c_4 \kappa^2 c_1^2+g_1 c_4 \kappa^3 c_1^3+6 g_2 \kappa c_1 c_4+3 g_2 \kappa^2 c_1^2 c_2^3-4 g_2 c_4 \kappa^2 c_1^2+g_2 c_4 \kappa^3 c_1^3) e_n^4+\cdots +O(e_n^7)
    \end{align}
Clearly, if $g_0+g_1+g_2=1$ then the coefficients of $e_n^2$ and $e_n^3$ are zero and hence the convergence order is four which is optimal in the 
sense of Kung and Traub. 
\end{proof}

\begin{remark}
 In (\ref{eq19}), if we  equate $g_0=0$, $g_1=0$, $g_2=1$ and define $G_2 (t_1,t_2) = 1+t_1+t_2+ \alpha t_1^2 + \beta t_2^2 $, then the resulted iterative scheme is the first two 
 steps of iterative method (\ref{eq9}). Similarly if we equate $g_0=0$, $g_2=0$, $g_1=1$, $\kappa = -1$ and define $G_1(t_2) = (1- {\lambda}^{-1} t_2)^{-\lambda}$, then 
 it reproduces the first three step of iterative schemes (\ref{eq10}) and (\ref{eq11}).
\end{remark}
In order to extend the iterative scheme (\ref{eq19}), we add one more step and get (\textbf{FD4})
\begin{align} \label{eq42}
 \begin{cases}
    w_n     & = x_n-\kappa f(x_n), \\
  y_n     & = x_n - \frac{f(x_n)}{f[x_n,w_n]},\\
  z_n     & = y_n -  \Big[ \frac{(1-g_1-g_2)}{f[y_n,w_n]}  G_0(t_1)+\frac{g_1}{f[y_n,x_n]}  G_1(t_2)+ \frac{g_2}{f[x_n,w_n]} G_2(t_1,t_2) \Big],\\
 x_{n+1} &= z_n- \Big[ \frac{(1-h_1-h_2-h_3-h_4-h_5)}{f[y_n,z_n]} S_0(t_3,t_4,t_5)+ \frac{h_1}{f[z_n,w_n]} S_1(t_1,t_3,t_4,t_5) 
 + \frac{h_2}{f[z_n,x_n]}
       S_2(t_2,t_3,t_4,t_5) \\
       &+ \frac{h_3}{f[x_n,y_n]} S_3 (t_2,t_3,t_4,t_5)  + \frac{h_4}{f[y_n,w_n]} S_4 (t_1,t_3,t_4,t_5)+ \frac{h_5}{f[x_n,w_n]} S_5(t_1,t_2,t_3,t_4,t_5) \Big] f(z_n),
 \end{cases}
\end{align}
where $t_1=\frac{f(y_n)}{f(x_n)}$, $t_2=\frac{f(y_n)}{f(w_n)}$,  $t_3=\frac{f(z_n)}{f(x_n)}$,  $t_4=\frac{f(z_n)}{f(w_n)}$, 
$t_5 =\frac{f(z_n)}{f(y_n)}$ and $h_1$, $h_2$, $h_3$, $h_4$, $h_5$, $\kappa\neq 0 $  are real numbers.   
The convergence of (\ref{eq42}) is discussed in the following theorem.
\newtheorem{mydef4}[mydef1]{Theorem}
\begin{mydef4}
 By including the statement of (Theorem 3), iterative scheme  (\ref{eq42}) has convergence order at least six if 
 \begin{align} \label{eq43}
 \begin{cases}
 & S_0(0,0,0)=1,\ S_1(0,0,0,0)=1,\ S_2(0,0,0,0)=1,\ S_3(0,0,0,0)=1,\ S_4(0,0,0,0)=1,\\
 & S_5(0,0,0,0,0)=1,\ \partial S_1 / \partial t_1 (0,0,0,0) = 1,\ \partial S_2 / \partial t_2 (0,0,0,0) = 1,\\
 &\partial S_3 / \partial t_2 (0,0,0,0) = 1,\ \partial S_4 / \partial t_1 (0,0,0,0) = 1,\\
  &\partial S_5 / \partial t_1 (0,0,0,0,0) = 1,\ \partial S_5 / \partial t_2 (0,0,0,0,0) = 1 
  \end{cases}
 \end{align}
 and all higher order derivatives of all functions  are bounded.
 \end{mydef4}
 
 \begin{proof}
   Taylor's expansion of various orders for $S_i$, $i=0,1,\cdots,5$ is given below.
  \begin{align}\label{eq44} 
  \begin{cases}
   S_0(t_3,t_4,t_5) &=  1+a_1 t_3 + a_2 t_4 +a_3 t_5+ \text{higher powers in $t_3$, $t_4$, $t_5$},  \\
   S_1(t_1,t_3,t_4,t_5) &= 1+t_1 + b_1 t_1^2+ b_2t_3+ b_3t_4+b_4t_5+ \text{higher powers in $t_1$, $t_3$, $t_4$, $t_5$},\\
   S_2(t_1,t_3,t_4,t_5) &= 1+t_2 + b_5 t_2^2+ b_6t_3+ b_7t_4+b_8t_5+ \text{higher powers in  $t_2$, $t_3$, $t_4$, $t_5$}, \\
   S_3(t_2,t_3,t_4,t_5) &= 1+t_2 + b_9 t_2^2+ b_{10}t_3+ b_{11}t_4+b_{12}t_5+ \text{higher powers in $t_2$, $t_3$, $t_4$, $t_5$}, \\
   S_4(t_1,t_3,t_4,t_5) &= 1+t_1 + b_{13} t_1^2+ b_{14}t_3+ b_{15}t_4+b_{16}t_5+ \text{higher powers in $t_1$, $t_3$, $t_4$, $t_5$}, \\
   S_5(t_1,t_2,t_3,t_4,t_5) &= 1+ t_1 + t_2 + b_{17} t_1^2 + b_{18} t_2^2+ b_{19}t_3+b_{20}t_4+b_{21}t_5
    + \text{higher powers in $t_1$, $t_2$, $t_3$, $t_4$, $t_5$},
  \end{cases}
  \end{align}
  where $a_1$, $a_2$, $a_3$, $b_{1}$, $b_{2}$, $b_3$, $b_4$, $b_5$, $b_6$, $b_7$, $b_8$, $b_9$, $b_{10}$, $b_{11}$, $b_{12}$, $b_{13}$, $b_{14}$, $b_{15}$, $b_{16}$, $b_{17}$, $b_{18}$, $b_{19}$, $b_{20}$,  and $b_{21}$ are real numbers. 
  By using calculations  from (\ref{eq35})-(\ref{eq41}) for $ x_{n+1} = z_n- \Big[ \frac{(1-h_1-h_2-h_3-h_4-h_5)}{f[y_n,z_n]} S_0(t_3,t_4,t_5)+ \frac{h_1}{f[z_n,w_n]} S_1(t_1,t_3,t_4,t_5)  \\
  + \frac{h_2}{f[z_n,x_n]} S_2(t_2,t_3,t_4,t_5)+ \frac{h_3}{f[x_n,y_n]} S_3 (t_2,t_3,t_4,t_5)  + \frac{h_4}{f[y_n,w_n]} S_4 (t_1,t_3,t_4,t_5) 
       + \frac{h_5}{f[x_n,w_n]} S_5(t_1,t_2,t_3,t_4,t_5) \Big] f(z_n)$, we obtain the following error equation.
 \begin{align} \label{eq45}
        e_{n+1} &= -(2 h_4 b_{13} c_2^2 \kappa c_1-h_4 b_{13} c_2^2 \kappa^2 c_1^2+ a_3 g_1 M_1 c_2^2-h_5 b_{18} c_2^2-h_5 b_{17} c_2^2-h_3 b_9 c_2^2+h_5 b_{21} g_1 \kappa^2 c_1^2 c_2^2\nonumber \\
        &+h_5 b_{21} g_2 c_2^2 \kappa c_1-2 h_5 b_{21} L_1 c_2^2 \kappa c_1+h_5 b_{21} L_1 c_2^2 \kappa^2 c_1^2+5 c_2^2 h_5+h_4 \kappa^2 c_1^2 c_2^2+h_4 \kappa c_1 c_3-4 h_4 c_2^2 \kappa c_1+h_1 \kappa^2 c_1^2 c_2^2\nonumber \\
        &+h_1 \kappa c_1 c_3-3 h_1 c_2^2 \kappa c_1-c_2^2 h_2 \kappa c_1-2 c_2^2 h_3 \kappa c_1-5 c_2^2 h_5 \kappa c_1+h_2 \kappa c_1 c_3+h_3 \kappa c_1 c_3+h_5 \kappa^2 c_1^2 c_2^2+h_5 \kappa c_1 c_3\nonumber \\
        &-a_3 \kappa^2 c_1^2 c_2^2- a_3 \kappa c_1 c_3+4  a_3 c_2^2 \kappa c_1+2  a_3 h_1 g_2 c_2^2+ \cdots +h_3 b_{12} L_1 c_2^2 \kappa^2 c_1^2-2 h_2 b_8 g_1 c_2^2 \kappa c_1+h_2 b_8 g_1 \kappa^2 c_1^2 c_2^2\nonumber \\
        &+h_2 b_8 g_2 c_2^2 \kappa c_1-2 h_2 b_8 L_1 c_2^2 \kappa c_1+h_2 b_8 L_1 c_2^2 \kappa^2 c_1^2) c_2 (-1+\kappa c_1) (-g_1 \kappa^2 c_1^2 c_2^2-L_1 c_2^2 \kappa^2 c_1^2-g_2 N_1 c_2^2 \kappa^2 c_1^2\nonumber \\
        &+\kappa^2 c_1^2 c_2^2+g_1 L_1 c_2^2 \kappa^2 c_1^2+g_2 L_1 c_2^2 \kappa^2 c_1^2-4 c_2^2 \kappa c_1+2 g_1 c_2^2 \kappa c_1+2 g_2 N_1 c_2^2 \kappa c_1-g_2 c_2^2 \kappa c_1-2 g_1 L_1 c_2^2 \kappa c_1\nonumber \\
        &+\kappa c_1 c_3+2 L_1 c_2^2 \kappa c_1-2 g_2 L_1 c_2^2 \kappa c_1-c_3-g_2 N_2 c_2^2-g_2 N_1 c_2^2+2 g_2 c_2^2+3 c_2^2-g_1 M_1 c_2^2-L_1 c_2^2+g_1 L_1 c_2^2\nonumber \\
        &+g_2 L_1 c_2^2) e_n^6+O(e_n^7) 
\end{align}
\end{proof}

\begin{remark}
 In (\ref{eq42}), by  equating $\kappa=-1$, $g_1=1$, $g_2=0$, $h_3=1$, $h_1=h_2=h_4=h_5=0$ and defining $G_2(t_2) = (1-\lambda^{-1}t_2)^{-\lambda}$, $S_3 (t_2,t_3,t_4,t_5)= (1-\lambda^{-1}t_2)^{-\lambda}$ 
 gives the iterative scheme (\ref{eq10}), by setting $\kappa=-1$, $g_1=1$, $g_2=0$, $h_1=h_2=h_3=h_4=h_5=0$ and defining $G_2(t_2) = (1-\lambda^{-1}t_2)^{-\lambda}$,$S_0(t_3,t_4,t_5)= 1$, 
 we obtain the iterative schemes (\ref{eq11}). (\ref{eq9}) can be constructed by choosing parameters $g_1=0$, $g_2=1$, $h_5=1$, $h_1=h_2=h_3=h_4=0$
 and weight functions $G_2(t_1,t_2)=1+t_1+t_2+\alpha t_1^2 + \beta t_2^2 $, $S_5(t_1,t_2,t_3,t_4,t_5)=1+t_1+t_2+\alpha t_1^2 + \beta t_2^2 + \eta t_5$. 
\end{remark}
 F. Soleymani constructs derivative-free sixth-order iterative scheme  (\ref{eq14}). In the following theorem, we present a modified iterative
 family with convergence order seven.
\newtheorem{mydef5}[mydef1]{Theorem}
\begin{mydef5}
 Let  $f: D\subseteq \Re \rightarrow \Re $ be a sufficiently differentiable function, and $\alpha \in D$ is a simple root of $f(x)=0$, 
 for an open interval $D$. If $x_0$ is chosen sufficiently close to $\alpha$, then iterative scheme  (\textbf{FD5})
\begin{align}  \label{eq46}
 \begin{cases}
  w_n &= x_n - \kappa f(x_n),  \\
  y_n &= x_n - \frac{f(x_n)}{f[x_n,w_n]}, \\
  z_n &= y_n - \Bigg[ \frac{1-g_1-g_2}{f[y_n,w_n]} G_0(t_1) + \frac{g_1}{f[x_n,y_n]} G_1(t_2) + \frac{g_2}{f[x_n,w_n]} G_2(t_1,t_2) \Bigg] f(y_n), \\
  x_{n+1} &= z_n - H(t_3,t_4,t_5) \frac{f(z_n)}{f[z_n,y_n]-(h-1)f[z_n,w_n]+h f[z_n,x_n]-h f[y_n,x_n] +(h-1) f[y_n,w_n]}, 
 \end{cases}
\end{align}
where
\begin{align}\label{eq47}
 \begin{cases}
  G_0(t_1) & = 1+t_1 + a_1 t_1^2+ \text{higher powers in $t_1$}, \\
  G_1(t_2) & = 1+t_2 + a_2 t_2^2+ \text{higher powers in $t_2$}, \\
  G_2(t_1,t_2) & = 1+t_1+t_2 + a_3 t_1^2 + a_4 t_2^2 + a_5 t_1 t_2+ \text{higher powers in $t_1$, $t_2$},  \\
  H(t_3,t_4,t_5) & = 1+a_6 t_3 + a_7 t_4 + a_8 t_5+ \text{higher powers in $t_3$, $t_4$, $t_5$}, 
 \end{cases}
\end{align}
$t_i$, $i=1,2,..,5$ are consistent with previous definition in (\ref{eq42}), $g_1$, $g_2$, $h$, $a_j$, $j=1,\cdots,8$ are real numbers. 
\end{mydef5}
has convergence order at-least six if 
\begin{align} \label{eq48}
 \begin{cases}
  G_0(0)=1,\ G_1(0)=1,\ G_2(0,0)=1,\ H(0,0,0)=1, \\
  dG_0/dt_1(0) = 1,\ dG_1/dt_2(0) = 1,\ \partial G_2/\partial t_1(0,0) = 1,\ \partial G_2/\partial t_2(0,0) = 1,
 \end{cases}
\end{align}
and has convergence order seven if along with (\ref{eq48}), the following condition is valid
\begin{align}\label{eq49}
 \partial H/\partial t_5 (0,0,0) = 0.
\end{align}

\begin{proof}
The terminology of (Theorem 3)  and  expressions in  (\ref{eq35}), (\ref{eq36}) assist to perform  all necessary
calculations by using Maple 13 \cite{23} software-package in order to get the following error equation:
\begin{align}\label{eq50}
 e_{n+1} &=  \boldsymbol{a_8} c_2    (-1+\kappa c_1) (\kappa^2 c_1^2 c_2^2-c_3+\kappa c_1 c_3+3 c_2^2-4 c_2^2 \kappa c_1+g_1 a_1 c_2^2-a_1 c_2^2+2 a_1 c_2^2 \kappa c_1-a_1 c_2^2 \kappa^2 c_1^2-g_2 c_2^2 \kappa c_1\nonumber \\
 & +2 g_1 c_2^2 \kappa c_1-g_1 \kappa^2 c_1^2 c_2^2+g_2 a_1 c_2^2-g_1 a_2 c_2^2-g_2 a_3 c_2^2-g_2 a_5 c_2^2-g_2 a_4 c_2^2+2 g_2 c_2^2-2 g_1 a_1 c_2^2 \kappa c_1+g_1 a_1 c_2^2 \kappa^2 c_1^2\nonumber \\
 &-2 g_2 a_1 c_2^2 \kappa c_1+g_2 a_1 c_2^2 \kappa^2 c_1^2+2 g_2 a_3 c_2^2 \kappa c_1-g_2 a_3 c_2^2 \kappa^2 c_1^2+g_2 a_5 c_2^2 \kappa c_1)^2 \boldsymbol{e_n^6} \nonumber \\
 & +  -(\kappa^2 c_1^2 c_2^2-c_3+\kappa c_1 c_3+3 c_2^2-4 c_2^2 \kappa c_1+g_1 a_1 c_2^2-a_1 c_2^2+2 a_1 c_2^2 \kappa c_1-a_1 c_2^2 \kappa^2 c_1^2-g_2 c_2^2 \kappa c_1+2 g_1 c_2^2 \kappa c_1\nonumber \\
 & -g_1 \kappa^2 c_1^2 c_2^2+g_2 a_1 c_2^2-g_1 a_2 c_2^2-g_2 a_3 c_2^2-g_2 a_5 c_2^2-g_2 a_4 c_2^2+2 g_2 c_2^2-2 g_1 a_1 c_2^2 \kappa c_1+g_1 a_1 c_2^2 \kappa^2 c_1^2\nonumber \\
 &+\cdots+\nonumber \\
 &-a_6 a_1 c_2^4-2 a_8 c_3^2-30 a_8 c_2^4+3 a_6 c_2^4-3 a_6 \kappa^2 c_1^2 c_2^2 c_3+a_6 c_3 \kappa^3 c_1^3 c_2^2+3 a_6 \kappa c_1 c_2^2 c_3-5 a_6 g_1 \kappa^2 c_1^2 c_2^4\nonumber \\
 &+2 a_6 g_1 \kappa c_1 c_2^4+4 a_6 g_1 \kappa^3 c_1^3 c_2^4+4 a_6 g_2 \kappa^2 c_1^2 c_2^4-5 a_6 g_2 \kappa c_1 c_2^4-a_6 g_2 \kappa^3 c_1^3 c_2^4-a_6 g_1 \kappa^4 c_1^4 c_2^4-a_6 a_1 c_2^4 \kappa^4 c_1^4\nonumber \\
 &+4 a_6 a_1 c_2^4 \kappa^3 c_1^3)\boldsymbol{e_n^7} +O(e_n^8) 
\end{align}
Clearly, $a_8=\partial H/\partial t_5 (0,0,0)$ is a factor in asymptotic error constant of $e_n^6$, if $a_8=0$ then (\ref{eq47}) has 
convergence order  seven otherwise convergence order is six.
\end{proof}
\begin{remark}
 Iterative scheme (\ref{eq47})  reproduces F. Soleymani derivative-free sixth-order iterative scheme (\ref{eq14}) if 
 \begin{align}\label{eq51}
  h&=0,\ g_1=g_2=0,\ G_0(t_1)=1,\ H(t_3,t_4,t_5)=1. 
 \end{align}
\end{remark}
The derivative-free version of (\ref{eq17}) is (\textbf{FD6}):
\begin{align}\label{eq52}
 \begin{cases}
  w_n &= x_n- \kappa f(x_n), \\
  y_n &= x_n - \frac{f(x_n)}{f[x_n,w_n]},\\
  z_n &= y_n - \Bigg[ \frac{1-g_1-g_2}{f[y,w]} G_0(t_1) + \frac{g_1}{f[x_n,y_n]} G_1(t_2) + \frac{g_2}{f[x_n,w_n]} G_2(t_1,t_2)\Bigg],\\
  x_{n+1} &= z_n - H(t_3,t_4,t_5) \frac{f(z_n)}{f[z_n,y_n]+ (f[z_n,x_n]-f[w_n,x_n])\frac{z_n-y_n}{z_n-x_n}},
 \end{cases}
\end{align}
where
\begin{align} \label{eq53}
 \begin{cases}
  G_0(t_1) & = 1+t_1 + a_1 t_1^2+ \text{higher powers in $t_1$}, \\
  G_1(t_2) & = 1+t_2 + a_2 t_2^2+ \text{higher powers in $t_2$}, \\
  G_2(t_1,t_2) & = 1+t_1+t_2 + a_3 t_1^2 + a_4 t_2^2 + a_5 t_1 t_2+ \text{higher powers in $t_1$, $t_2$},  \\
  H(t_3,t_4,t_5) & = 1+a_6 t_3 + a_7 t_4 + a_8 t_5+ \text{higher powers in $t_3$, $t_4$, $t_5$}, 
 \end{cases}
\end{align}
$t_i$, $i=1,2,\cdots,5$ are consistent with previous definition in (\ref{eq42}), $g_1$, $g_2$, $h$, $a_j$, $j=1,\cdots,8$ are real numbers. 
\newtheorem{mydef6}[mydef1]{Theorem}
\begin{mydef6}
 Let  $f: D\subseteq \Re \rightarrow \Re $ be a sufficiently differentiable function, and $\alpha \in D$ is a simple root of $f(x)=0$, 
 for an open interval $D$. If $x_0$ is chosen sufficiently close to $\alpha$, then iterative scheme (\ref{eq53})
has convergence order at-least six if  
\begin{align} \label{eq54}
 \begin{cases}
  G_0(0)=1,\ G_1(0)=1,\ G_2(0,0)=1,\ H(0,0,0)=1, \\
  dG_0/dt_1(0) = 1,\ dG_1/dt_2(0) = 1,\ \partial G_2/\partial t_1(0,0) = 1,\ \partial G_2/\partial t_2(0,0) = 1.
 \end{cases}
\end{align}
\end{mydef6}
\begin{proof}
We obtain the  following error equation
\begin{align}\label{eq55}
 e_{n+1} &=  (-1+\kappa c_1) (\kappa^2 c_1^2 c_2^2-c_3+\kappa c_1 c_3+3 c_2^2-4 c_2^2 \kappa c_1+g_1 a_1 c_2^2+g_2 a_1 c_2^2-g_1 a_2 c_2^2-g_2 c_2^2 a_5-g_2 a_4 c_2^2\nonumber \\
 &-g_2 a_3 c_2^2-a_1 c_2^2-2 g_1 a_1 c_2^2 \kappa c_1+g_1 a_1 c_2^2 \kappa^2 c_1^2-2 g_2 a_1 c_2^2 \kappa c_1+g_2 a_1 c_2^2 \kappa^2 c_1^2+2 g_2 a_3 c_2^2 \kappa c_1+2 a_1 c_2^2 \kappa c_1\nonumber \\
 &-a_1 c_2^2 \kappa^2 c_1^2+2 g_1 c_2^2 \kappa c_1-g_1 \kappa^2 c_1^2 c_2^2-g_2 c_2^2 \kappa c_1+2 g_2 c_2^2-g_2 a_3 c_2^2 \kappa^2 c_1^2+g_2 c_2^2 a_5 \kappa c_1) (\kappa^2 c_1^2 c_2^2-c_2^2 \kappa c_1+3 c_2^2 a_8\nonumber \\
 &-a_8 c_3+2 c_2^2 a_8 g_2-a_1 c_2^2 a_8-c_2^2 a_8 a_2 g_1-c_2^2 a_8 g_2 a_3-c_2^2 a_8 g_2 a_5-c_2^2 a_8 g_2 a_4+c_1 a_8 c_3 \kappa+a_1 c_2^2 a_8 g_1\nonumber \\
 &+a_1 c_2^2 a_8 g_2-c_1^2 c_2^2 a_8 \kappa^2 g_1+2 c_1 c_2^2 a_8 \kappa g_1-c_1 c_2^2 \kappa a_8 g_2-c_1^2 c_2^2 \kappa^2 a_8 g_2 a_3+c_1 c_2^2 \kappa a_8 g_2 a_5+2 c_1 c_2^2 \kappa a_8 g_2 a_3\nonumber \\
 &-a_1 c_1^2 c_2^2 a_8 \kappa^2+2 a_1 c_1 c_2^2 a_8 \kappa+a_1 c_1^2 c_2^2 a_8 \kappa^2 g_1+a_1 c_1^2 c_2^2 \kappa^2 a_8 g_2-2 a_1 c_1 c_2^2 a_8 \kappa g_1-2 a_1 c_1 c_2^2 \kappa a_8 g_2+c_1^2 c_2^2 a_8 \kappa^2\nonumber \\
 &-4 c_1 c_2^2 a_8 \kappa) c_2  e_n^6 +O(e_n^7)
\end{align}
\end{proof}

\section{Numerical computations}

\newtheorem{mydef7}{Definition}
\begin{mydef7}
 Let $x_{n-1}$, $x_n$ and $x_{n+1}$ be successive iterations closer to the root $\alpha$ of $f(x)=0$, the computational order of convergence (COC) \cite{24}, can be approximated by
 \begin{align} \label{eq56}
 COC &\approx  \frac{ln|(x_{n+1}-\alpha)(x_n-\alpha)^{-1}|}{ln|(x_n-\alpha)(x_{n-1}-\alpha)^{-1}|}.
 \end{align}
\end{mydef7}

A set of twelve functions is listed in Table 1 which  is taken from \cite{25} to validate the iterative methods and performance.
In all methods, twelve  total number of function evaluations (TNFE) are used.  

\subsection{Sixth-order convergence performance evaluation of \textbf{FD1} }
For the purpose of comparison between \textbf{FD1} and \textbf{SG}, \textbf{NT1}, \textbf{NT2}, \textbf{CH}, \textbf{GR}, \textbf{AL}, we derive  an iterative method 
from (\ref{eq15}) which is given as (\textbf{FD1-M1} )
\begin{align} \label{eq57}
 \begin{cases}
  y_n &= x_n - \frac{f(x_n)}{f'(x_n)}, \\
  z_n &= y_n -  G(t_1) \frac{f(y_n)}{f'(x_n)},\\
  x_{n+1} &= z_n -  H_1(t_1) H_2(t_2) H_3(t_3)  \frac{f(z_n)}{f'(x_n)},
 \end{cases}
\end{align}
where $t_1=\frac{f(y_n)}{f(x_n)}, t_2 = \frac{f(z_n)}{f(x_n)}, t_3 = \frac{f(z_n)}{f(y_n)}, G(t_1) =  \frac{1}{1-2t_1}, H_1(t_1) = \frac{1}{1-2t_1-t_1^2}
, H_2(t_2) = 1+2t_2$ and $H_3(t_3) = \frac{1}{1-1.1t_3}$. Table 2 shows the performance of different listed iterative methods 
in terms of absolute error ($|x_n-\alpha|$). Clearly \textbf{FD1-M1} is competitive with other developed methods and freedom to choose 
different weight functions, makes it more accurate and stable in the third step of iterative scheme (\ref{eq15}).
\begin{table}[!htbp]
\caption{List of test functions }\label{tb1: test functions }
\begin{center}
 \begin{tabular}{p{6.5cm}  p{2.75cm}}
  \hline  
  Functions & Roots \\
  \hline 
$f_1 (x) = exp(x)\ sin(x) + ln(1 + x^2 )$    		&  $\alpha=0 $ \\ 
$f_2 (x) = x^{15} + x^4 + 4x^2 -15  		$	&  $\alpha = 1.148538 . . .$  \\
$f_3 (x) = (x - 2)(x^{10} + x + 1)\ exp(-x-1)$		&  $\alpha=2  $\\
$f_4 (x) = exp(-x^2 + x + 2) - cos(x + 1) + x^3 + 1$	&  $\alpha = -1$  \\
$f_5 (x) = (x+1)\ exp(sin(x))-x^2 \ exp(cos(x))-1$	&  $\alpha=0  $\\
$f_6 (x) = sin(x)^2 - x^2 + 1		$		&  $\alpha = 1.40449165 . . . $ \\
$f_7 (x) = 10\ exp(-x^2 ) -1		$		&  $\alpha = 1.517427 . . .  $\\
$f_8 (x) = (x^2 -1)^{-1} - 1		$		&  $\alpha = 1.414214 . . .  $\\
$f_{9} (x) = ln(x^2 + x + 2) - x + 1	$		&  $\alpha = 4.15259074 . . .$  \\
$f_{10} (x) = cos(x)^2 - x/5		$	        &  $\alpha = 1.08598268 . . .$  \\
$f_{11} (x) = x^{10} - 2x^3 - x + 1	$       	&  $\alpha = 0.591448093 . . .$  \\
$f_{12} (x) = exp(sin(x)) - x + 1      $               &  $\alpha = 2.63066415 . . . $\\
\hline
\end{tabular}
\end{center}
\end{table}

\begin{table}[!htbp]
\caption{Numerical comparison between absolute errors($|x_n-\alpha|,TNFE=12) $ }\label{tb2: Numerical comparison.}
\begin{center}
 \begin{tabular}{p{1.3cm}  p{1.6cm}  p{1.6cm} p{1.6cm} p{1.6cm}  p{1.6cm} p{1.6cm}  p{1.6cm} }
  \hline  
   $f_n(x)$,$x_0$  & FD1-M1 & SG & NT1 & NT2 & CH & GR & AL  \\  
        &   &  $a=-1$ & $a=-1$ &  &  &  & $\theta=-1.01$  \\
  \hline 
    $f_1$,\ 0.25    &   \textbf{2.71 e-142}     & 2.82e-121   &  1.22e-127  &  4.25e-84     &  5.12e-82    &    9.23 e-93    &   6.22e-141     \\ 
   $ f_2$ ,\ 1.1    &   \textbf{3.90 e-190}  & 3.72 e-153   & 1.23 e-122  &  1.47 e-95   &  2.45 e-104    &  4.43 e-141   &  1.60 e-123    \\ 
   $f_3$,\ 2.1      &   \textbf{2.31e-181}   &  3.75e-135   &  6.27e-135   &   2.89e-93    &   1.07e-94    &   6.52e-109    &   1.80e-147    \\ 
   $f_4 $,\ -0.5    &   \textbf{5.28 e-218}   &  8.50 e-181   &  5.30 e-143   &   2.70 e-173    &   1.27 e-174    &   5.89 e-207    &   8.99 e-191   \\ 
   $f_5 $,\ 0.25    &   \textbf{1.96e-247}   &  4.05e-184   &  2.51e-165   &   3.77e-213    &   7.33e-162    &   2.97e-243    &   2.20e-215     \\ 
   $f_6 $,\ 1.2     &   \textbf{3.79e-215}    &  1.41e-181    &  7.22e-173    &   3.90e-123     &   1.69e-123     &   6.80e-154     &   2.42e-180    \\ 
   $f_7 $,\ 1.0     &  \textbf{5.26e-123}   & 2.88e-95   & 2.10e-113   &  3.87e-69    &  6.28e-67    &  1.81e-76    &  2.55e-118    \\ 
   $f_8 $,\ 1.6     &  \textbf{2.27e-181}   & 1.18e-157   & 1.65e-97   &  6.69e-59    &  1.23e-103    &  2.26e-147    &  4.14e-60    \\ 
   $f_{9} $,\ 4.4   &  2.92e-402  &  3.13e-366  &  4.80e-387  &   2.35e-344   &   6.52e-341   &   1.89e-353   &  \textbf{ 2.41e-432 }  \\ 
   $f_{10} $,\ 1.5  &  2.33e-56  &  2.22e-53  &  3.42e-64  &   3.58e-31   &   10.1   &   5.45e-56   &  \textbf{ 8.42e-58 }  \\ 
   $f_{11}$,\ 0.25  &  \textbf{8.99e-115}  &  1.30e-97  &  1.54e-91  &  1.31e-53   &   3.86e-45   &   4.77e-89   &   1.72e-92   \\ 
   $f_{12}$,\ 2.0   & \textbf{3.53e-173}   &   1.80e-132   &   3.42e-138  &    1.84e-117    &    2.78e-115    &    3.52e-126    &    8.05e-162  \\ 
\hline
\end{tabular}
\end{center}
\end{table}

\begin{table}[!htbp]
\caption{Computational order of convergence (COC)  }\label{tb3: COC.}
\begin{center}
 \begin{tabular}{p{1.5cm}  p{1.5cm}  p{1.5cm} p{1.5cm} p{1.5cm}  p{1.5cm} p{1.5cm}  p{1.5cm} }
  \hline  
  $f_n(x)$ & FD1-M1 & SG & NT1 & NT2 & CH & GR & AL  \\  
  \hline 
   $f_1$        &  5.9999 & 5.9994 &  6.0003 & 5.9989 & 5.9987 & 5.9995 & 5.9965  \\ 
   $f_2$        &  6.0000 & 6.0000 &  5.9999 & 5.9997 & 5.9998 & 6.0000 & 5.9933  \\ 
   $f_3$        &  6.0000 & 6.0000 &  5.9999 & 5.9995 & 5.9996 & 5.9999 & 6.0012  \\ 
   $f_4 $       &  6.0000 & 6.0000 &  6.0000 & 6.0000 & 6.0000 & 6.0000 & 5.9998  \\ 
   $f_5 $       &  6.0000 & 6.0000 &  5.9999 & 6.0000 & 5.9999 & 6.0000 & 5.9999  \\ 
   $f_6 $       &  6.0000 & 6.0000 &  5.9999 & 5.9999 & 5.9999 & 6.0000 & 5.9998  \\ 
   $f_7 $       &  5.9999 & 5.9992 &  6.0001 & 5.9974 & 5.9972 & 5.9991 & 5.9860  \\ 
   $f_8 $       &  6.0000 & 6.0000 &  6.0001 & 5.9985 & 5.9999 & 6.0000 & 6.0971  \\ 
   $f_{9} $     &  6.0000 & 6.0000 &  6.0000 & 6.0000 & 6.0000 & 6.0000 & 6.0000  \\ 
   $f_{10} $    &  5.9806 & 6.0038 &  6.0098 & 5.8343 & 3.4125 & 5.9880 & 6.1207  \\ 
   $f_{11} $    &  5.9999 & 5.9981 &  6.0001 & 5.9909 & 5.9856 & 6.0000 & 6.0120  \\ 
   $f_{12}$     &  6.0000 & 6.0000 &  6.0000 & 6.0002 & 6.0002 & 6.0001 & 5.9981  \\
\hline
\end{tabular}
\end{center}
\end{table}

\subsection{Seventh-order convergence performance evaluation of \textbf{FD1} }
To compare the results, we construct seventh-order iterative scheme (\textbf{FD1-M2}) from (\ref{eq15}):
\begin{align} \label{eq58}
 \begin{cases}
  y_n &= x_n - \frac{f(x_n)}{f'(x_n)}, \\
  z_n &= y_n -  G(t_1) \frac{f(y_n)}{f'(x_n)},\\
  x_{n+1} &= z_n -  H_1(t_1) H_2(t_2) H_3(t_3)  \frac{f(z_n)}{f'(x_n)},
 \end{cases}
\end{align}
where $t_i's$ are defined previously, $G(t_1) =  \frac{1}{1-2t_1}, H_1(t_1) = \frac{1}{1-2t_1-t_1^2}
, H_2(t_2) = 1+2.1t_2$ and $H_3(t_3) = \frac{1}{1-t_3}$. In (\ref{eq17}), if we define
$A(t_1)=\frac{1+t_1}{1-t_1}$ then (\ref{eq17}) is a seventh-order iterative scheme (\textbf{AL1}) which is
developed in \cite{17}  and if 
$A(t_1)=\frac{1}{(1-t_1)^2}$ then again (\ref{eq17}) is a seventh-order iterative scheme (\textbf{FD2-M1}).

\begin{table}[!htbp]
\caption{Numerical comparison between absolute errors($|x_n-\alpha|,TNFE=12)$ and COC }\label{tb4: Numerical comparison.}
\begin{center}
 \begin{tabular}{p{1.3cm}  p{1.6cm}  p{1.6cm} p{1.6cm} p{1.6cm}  p{1.6cm} p{1.6cm}  p{1.6cm} }
  \hline  
   $f_n(x)$,$x_0$  & FD1-M2 &   FD2-M1 & AL1 & FD1-M2    & FD2-M1  & AL1   \\  
      &      &      &      &  (COC)   &  (COC) &  (COC)  \\  
  \hline 
   $f_1$,\ 0.25    &   \textbf{3.60e-182}  &   8.55e-177  &   6.80e-167  &   6.9998  &   6.9999  &   6.9997   \\ 
   $f_2$,\ 1.1     &   \textbf{2.99e-263}   &   6.71e-172   &   1.63e-145   &   7.0000   &   7.0000   &   6.9999   \\ 
   $f_3$,\ 2.1     &   \textbf{8.60e-223}  &  3.95e-177  &  3.11e-178  &  7.0000  &  7.0000  &  7.0000   \\ 
   $f_4 $,\ -0.5   &   \textbf{2.78e-297}   &  4.95e-247   &  1.45e-233   &  7.0000   &  7.0000   &  7.0000   \\ 
   $f_5 $,\ 0.25   &   \textbf{1.31e-307}   &  1.45e-306   &  2.81e-273   &  7.0000   &  7.0000   &  7.0000   \\ 
   $f_6 $,\ 1.2    &   \textbf{5.02e-304}  &  2.70e-246  &  1.27e-222  &  7.0000  &  7.0000  &  7.0000    \\ 
   $f_7 $,\ 1.0    &   \textbf{1.51e-156}  &  1.83e-145  &  3.98e-136  &  6.9997  &  6.9998  &  6.9996   \\ 
   $f_8 $,\ 1.6    &   \textbf{4.95e-240}  &  3.58e-102  &  2.35e-67  &  6.9995  &  7.0001  &  6.9978   \\ 
   $f_{9} $,\ 4.4  &   \textbf{5.10e-662}  &  3.74e-567  &  1.07e-561  &  7.0000  &  7.0000  &  7.0000  \\ 
   $f_{10} $,\ 1.5 &   1.08e-77  &  \textbf{1.33e-82}  &  1.73e-64  &  6.8846  &  7.0083  &  6.9701   \\ 
   $f_{11}$,\ 0.25 &   \textbf{8.67e-162}  &  4.51e-141  &  1.31e-103  &  7.0000  &  7.0003  &  7.0019   \\ 
   $f_{12}$,\ 2.0  &   \textbf{1.30e-242}  &  4.29e-196  &  9.45e-193  &  7.0000  &  7.0000  &  6.9999   \\ 
\hline
\end{tabular}
\end{center}
\end{table}
Table 4 shows absolute error due to method \textbf{FD2-M1} has small magnitude in comparison with other describe methods for the set of given
twelve test functions and initial guesses. Only in the case of $f_9$ and $f_{10}$, \textbf{AL} has better absolute error magnitude.

\subsection{Sixth-order convergence performance evaluation of \textbf{FD4, FD5, FD6} }
For comparison, we develop some particular cases of derivative-free methods. In all calculations the value of parameter $\kappa=1/100$ and 
$t_1=\frac{f(y_n)}{f(x_n)}$, $t_2=\frac{f(y_n)}{f(w_n)}$, $t_3=\frac{f(z_n)}{f(x_n)}$, $t_4=\frac{f(z_n)}{f(w_n)}$, $t_5=\frac{f(z_n)}{f(y_n)}$.
\begin{itemize}
 \item \textbf{FD4}
 \begin{align}\label{eq63}
  \begin{cases}
   w_n & = x_n -\kappa f(x_n), \\
   y_n &= x_n - \frac{f(x_n)}{f[x_n,w_n]},\\
   z_n &= y_n - \frac{1+t_1-t_2}{1-2t_2} \frac{f(y_n)}{f[x_n,w_n]},\\
   x_{n+1} &= z_n -\frac{1}{2}\Big[ \frac{1}{1-t_5-2t_3-2t_4} \frac{1}{f[z_n,y_n]} + \frac{1-2t_2}{1-3t_2} \frac{1}{f[z_n,x_n]}   \Big]f(z_n).
  \end{cases}
 \end{align}
 \item \textbf{FD5}
 \begin{align} \label{eq64}
  \begin{cases}
   w_n &= x_n-\kappa f(x_n), \\
  y_n &= x_n - \frac{f(x_n)}{f[x_n,w_n]},\\
  z_n &= y_n - \frac{1}{1-t_2} \frac{f(y_n)}{f[x_n,y_n]},\\
  x_{n+1} &= z_n - \Big( 1-\frac{t_5}{10}\Big)^{-1} \frac{f(z_n)}{f[z_n,y_n]+f[z_n,w_n]-f[y_n,w_n]}.
  \end{cases}
 \end{align}
\item  \textbf{FD6}
\begin{align} \label{eq65}
 \begin{cases}
  w_n &= x_n - \kappa f(x_n), \\
  y_n &= x_n - \frac{f(x_n)}{f[x_n,w_n]}, \\
  z_n &= y_n - \frac{1-t_1+t_2}{1-2t_1} \frac{f(y_n)}{f[x_n,w_n]},\\
  x_{n+1} &= z_n - \Big( 1-\frac{t_5}{10}\Big)^{-1}  \frac{f(z_n)}{f[z_n,y_n]+(f[z_n,x_n]-f[x_n,w_n]) \frac{z_n-y_n}{z_n-x_n}}. 
 \end{cases}
\end{align}
\item \textbf{SK2M1, SK2M2}\\
In \cite{18}, Authors developed two methods by selecting different values of parameters. If $\alpha=\beta=\frac{5}{2}$, $\eta= 1$ then \textbf{SK2} 
is \textbf{SK2M1} and if $\alpha=\beta=\eta=1$ then \textbf{SK2} is \textbf{SK2M2}.
\end{itemize}
Table 5 shows overall performance of derivative-free sixth-order iterative methods \textbf{FD4, FD5, FD6}. Absolute error $|x_n-\alpha|$, for
developed sixth-order derivative-free methods, is comparatively better than referenced sixth-order derivative-free iterative methods. 
Computational order of convergence is given in Table 6. \textbf{dgt} stands for \textbf{divergent} and \textbf{X} is for \textbf{no information}.

\begin{table}[!htbp]
\caption{ Numerical comparison between absolute errors($|x_n-\alpha|,TNFE=12)$  }\label{tb5: Numerical comparison.}
\begin{center}
 \begin{tabular}{p{1.3cm}  p{1.2cm}  p{1.2cm} p{1.2cm} p{1.2cm}  p{1.2cm} p{1.2cm}  p{1.2cm}  p{1.2cm}  p{1.2cm}  }
  \hline  
  $f_n(x)$ & FD4 & FD5 & FD6 & TS1 & TS2 & SK2M1 & SK2M2 & FS1 & FS2  \\  
  \hline 
    $f_1$, \ 0.25      & \textbf{4.e-172}&2.e-158&\textbf{8.e-172}&9.e-54&2.e-66&3.e-103&4.e-83&1.e-59&2.e-104	\\ 
   $f_2$, \  1.1       &1.e-190&2.e-297&6.e-221&0.01&3.e-14&7.e-162&7.e-241&2.e-7&\textbf{9.e-333}	\\ 
   $f_3$, \ 2.1        &2.e-107&\textbf{2.e-203}&4.e-168&3.e-18&4.e-24&2.e-197&6.e-145&3.e-7&\textbf{5.e-204}	\\ 
   $f_4 $, \-0.5       &\textbf{3.e-220}&2.e-198&8.e-213&3.e-127&2.e-202&1.e-217&3.e-175&9.e-196&4.e-171	\\ 
   $f_5 $, \ 0.25      &\textbf{3.e-247}&\textbf{7.e-248}&6.e-215&2.e-129&8.e-172&4.e-157&1.e-197&2.e-171&3.e-215	\\ 
   $f_6 $, \ 1.2       &1.e-200&2.e-186&\textbf{5.e-247}&1.e-132&5.e-153&8.e-112&2.e-161&8.e-194&2.e-161 	\\ 
   $f_7 $, \ 1.0       & 3.e-140&8.e-130&\textbf{8.e-146}&dgt &dgt&3.e-78&1.e-65&dgt&6.e-87	\\ 
   $f_8 $, \ 1.6       & \textbf{3.e-168}&1.e-90&\textbf{8.e-165}&5.e-55&2.e-71&2.&5.e-25&2.e-68&3.e-145	\\ 
   $f_{9} $, \ 4.4     &\textbf{5.e-466}&3.e-392&1.e-404&2.e-395&2.e-424&1.e-383&1.e-366&2.e-429&5.e-378 	\\ 
   $f_{10} $, \ 1.5    &1.e-70&6.e-70&8.e-65&1.e-138&\textbf{3.e-168}&0.06&1.&9.e-156&5.e-68 	\\ 
   $f_{11} $, \ 0.25   &7.e-118&4.e-118&\textbf{4.e-169}&2.e-71&1.e-105&1.&8.e-29&1.e-118&3.e-85  	\\ 
   $f_{12}$, \ 2.0     & 1.e-189&3.e-149&1.e-144&8.e-137&2.e-194&2.e-173&2.e-144&\textbf{5.e-215}&4.e-191 	\\
\hline
\end{tabular}
\end{center}
\end{table}

\begin{table}[!htbp]
\caption{Computational order of convergence (COC)  }\label{tb6: COC.}
\begin{center}
 \begin{tabular}{p{1.3cm}  p{1.2cm}  p{1.2cm} p{1.2cm} p{1.2cm}  p{1.2cm} p{1.2cm}  p{1.2cm}  p{1.2cm}  p{1.2cm}  }
  \hline  
  $f_n(x)$ & FD4 & FD5 & FD6 & TS1 & TS2 & SK2M1 & SK2M2 & FS1 & FS2  \\  
  \hline 
    $f_1$      &  6.01&6.00&6.0&5.99&6.00&5.98&6.00&5.98&6.00	\\ 
   $f_2$       &6.00&6.00&6.0&-1.15&6.24&6.00&6.00&3.33&6.00	\\ 
   $f_3$        &6.00&6.00&6.0&5.32&5.63&6.00&6.00&3.81&6.00	\\ 
   $f_4 $       &6.00&6.00&6.0&6.00&6.00&6.00&6.00&6.00&6.00	\\ 
   $f_5 $      &6.00&6.00&6.0&6.00&6.00&6.00&6.00&6.00&6.00	\\ 
   $f_6 $      & 6.00&6.00&6.0&6.00&6.00&6.00&6.00&6.00&6.00	\\ 
   $f_7 $       & 6.00&6.00&6.0&X&X&5.98&6.00&X&6.00	\\ 
   $f_8 $       & 6.01&6.00&6.0&5.99&6.00&5.23&5.78&6.00&6.00	\\ 
   $f_{9} $     & 6.00&6.00&6.0&6.00&6.00&6.00&6.00&6.00&6.00	\\ 
   $f_{10} $    &6.30&6.01&6.0&6.00&6.00&1.51&0.00299&6.00&5.99\\ 
   $f_{11} $   &  6.05&6.00&6.0&6.00&6.00&2.57&5.91&6.00&6.00		\\ 
   $f_{12}$     & 5.98&6.00&6.0&6.00&6.00&6.00&6.00&6.00&6.00\\
\hline
\end{tabular}
\end{center}
\end{table}

\begin{table}[!htbp]
\caption{Numerical comparison between absolute errors($|x_n-\alpha|,TNFE=12)$ and COC }\label{tb7: Numerical compariso}
\begin{center}
 \begin{tabular}{p{1.3cm}  p{2.7cm}  p{2.1cm} p{2.1cm} p{2.1cm}  p{2.1cm}  }
  \hline  
  $f_n(x)$ & FD7($\kappa$)(COC) & FS3-1(COC) & FS3-2(COC) & FS4-1(COC) & FS4-2(COC)   \\  
  \hline 
   $f_1$, \ 0.25       &  \textbf{2.45e-378}(1.0)(10)&3.23e-99(7)&4.83e-318(9)&2.96e-85(7)&6.60e-311(10)  \\ 
   $f_2$, \  1.1      &  \textbf{1.22e-388}(0.01)(7)&3.33e-31(6.42)&8.14e-24(7.02)&2.37e-7(3.30)&2.04(-0.0211) \\ 
   $f_3$, \ 2.1       &  \textbf{3.29e-266}(0.01)(7)&1.37e-40(6.80)&1.32e-26(7.42)&3.96e-8(4.18)&divergent(X)   \\ 
   $f_4 $, \-0.5       & \textbf{8.61e-262}(0.01)(7)&1.89e-221(7)&1.17e-144(7)&2.06e-210(7)&8.75e-75(6.97)   \\ 
   $f_5 $, \ 0.25      &  1.54e-303(0.01)(7)&2.47e-220(7)&\textbf{6.43e-397}(7)&1.33e-206(7)&1.25e-384(7)  \\ 
   $f_6 $, \ 1.2       &  \textbf{5.15e-353}(0.01)(7)&1.71e-251(7)&6.17e-74(7.01)&1.34e-257(7)&4.81e-44(7.05)   \\ 
   $f_7 $, \ 1.0       &  \textbf{3.65e-213}(0.01)(7)&0.320(4.96)&3.03(X)&dgt(X)&3.03(X)   \\ 
   $f_8 $, \ 1.6       & 1.10e-174(0.01)(7)&\textbf{1.09e-234}(7)&7.46e-22(7.09)&2.68e-233(7)&290(1.43)  \\ 
   $f_{9} $, \ 4.4     & \textbf{9.09e-689}(0.01)(7)&4.55e-617(7)&1.71e-517(7)&5.06e-640(7)&4.45e-506(7)  \\ 
   $f_{10} $, \ 1.5    & \textbf{3.81e-262}(-1.0)(7)&6.26e-174(7)&3.87e-5(2.04)&1.96e-214(7)&7.16(X)    \\ 
   $f_{11} $, \ 0.25   & \textbf{7.21e-212}(0.01)(7)&1.33e-125(7)&5.57e-77(6.99)&6.28e-147(7)&1.58e-52(6.97)  \\ 
   $f_{12}$, \ 2.0    &  \textbf{3.52e-271}(0.01)(7)&7.93e-257(7)&6.79e6(-7.79)&3.40e-222(7)&2.64(-0.261) \\
\hline
\end{tabular}
\end{center}
\end{table}

\subsection{Seventh-order convergence performance evaluation of \textbf{FD7 }}
\textbf{FD7} is seventh-order convergent iterative method which is deduced from  \textbf{FD5} and  is defined as
\begin{align}\label{eq66}
 \begin{cases}
  w_n &= x_n - \kappa f(x_n), \\
  y_n &= x_n - \frac{f(x_n)}{f[x_n,w_n]}, \\
  z_n &= y_n - \frac{1-2t_1}{1-3t_1} \frac{f(y_n)}{f[y_n,w_n]},\\
  x_{n+1} &= z_n - \frac{1}{1-t_3} \frac{f(z_n)}{f[z_n,y_n]+f[z_n,x_n]-f[x_n,y_n]}.
 \end{cases}
\end{align}
All  parameters set to zero value namely $\gamma=0$, $\delta=0$ in \textbf{FS3-1}, $\rho=0$, $\tau=0$ in \textbf{FS3-2} and $\omega=0$, 
$\phi=0$ in \textbf{FS4}. Table 7 represents the comparison between absolute errors as well as computational order of convergence. Clearly
\textbf{FD7} is superior than other methods in comparison. 

\section{Conclusions}
 We have presented  derivative-based and derivative-free iterative methods by introducing free weighting  parameters $g_0$, $g_1$, $g_2$ and 
 weight functions. Weighting parameters, in the second step of derivative-free iterative methods, give flexibility for selecting different forms 
 of approximation for first order derivative. A combination of derivative approximations can be constructed by fixing weighting parameters.
 The proper selection weight functions for both derivative-free and derivative-based iterative methods help to improve the accuracy. 
 The presented iterative methods are defined  in very general form which actually reproduce many existing iterative methods. 
 By introducing different form of weighting functions and selection of parameters produce new families of iterative methods.
 We explore some of them to show the effectiveness of newly constructed iterative methods.

\noindent\textbf{Acknowledgement.\ }
This research was supported by Spanish MICINN grants AYA2010-15685 and MEC grants AYA2010-15685, AYA2008-04211-C02-C01.

\end{document}